\documentclass{mcom-l}

\usepackage{graphicx}
\usepackage{graphics}
\usepackage{srcltx}
\usepackage{tikz}
\usepackage{float}
\usepackage{algorithm}
\usepackage{color}
\usepackage{algpseudocode}

\usepackage{amsmath}
\usepackage{amsthm}
\usepackage{amssymb}
\usepackage{colortbl}
\usepackage{epsfig}
\usepackage{amsfonts}
\usepackage{amsxtra}     

\usepackage{srcltx}
\usepackage{epsfig}
\usepackage{verbatim}

\usepackage{color}

\newtheorem{theorem}{Theorem}[section]
\newtheorem{lemma}[theorem]{Lemma}
\newtheorem{cor}[theorem]{Corollary}
\newtheorem{prop}[theorem]{Proposition}
\newtheorem{lem}[theorem]{Lemma}
\theoremstyle{definition}

\newtheorem{example}[theorem]{Example}

\newcommand{\tomlib}{\textsf{Tomlib}}
\newcommand{\gap}{\textsf{GAP}}
\theoremstyle{remark}

\numberwithin{equation}{section}

\begin{document}

\title{Computing The Table Of Marks Of A Cyclic Extension}


\author{L. Naughton}
\address{School of Mathematics, Statistics and Applied Mathematics, NUI, Galway}
\email{liam.naughton@nuigalway.ie, goetz.pfeiffer@nuigalway.ie}

\author{G. Pfeiffer}


\subjclass[2010]{Primary 20B40; Secondary 19A22, 20D30, 20D08, 20D10}

\date{May 2011}

\dedicatory{}

\begin{abstract}
The subgroup pattern of a finite groups $G$ is the table of marks of $G$ together with a list of 
representatives of the conjugacy classes of subgroups of $G$. In this article we present an algorithm for 
the computation of the subgroup pattern of a cyclic extension of $G$ from the subgroup pattern of $G$.  
Repeated application of this algorithm yields an algorithm for the
computation of the table of marks of a solvable group $G$, along a
composition series of~$G$.
\end{abstract}

\maketitle
\section{Introduction}
The actions of a finite group $G$ on finite sets $X$ are closely linked to the subgroup structure of $G$, since the isomorphism types of transitive 
actions of $G$ are in bijection to the conjugacy classes of subgroups of $G$. Thus properties of finite group actions have an impact on the subgroup 
structure of $G$, and vice versa. The correspondence between classes of subgroups of $G$ and transitive actions is made explicit in the \emph{table of marks} 
of $G$. This matrix was introduced by Burnside \cite{burn} as a tool to classify $G$-sets up to equivalence. In this context, 
the \emph{mark} 
of a subgroup $H$ of $G$ on $X$ is the number of fixed points
of $H$ in the action of $G$ on $X$, denoted by
$\beta_{X} (H)$.
If $H_{1}, \ldots H_{r}$ is a list of representatives of the subgroups of $G$ up to 
conjugacy, the table of marks of $G$ is then the $( r \times r)$-matrix
\[ \mathrm{M}(G) =(\beta_{G/H_{i}}(H_{j}))_{i, j=1,\ldots,r} .\]
Similar to the character table of $G$, which classifies matrix representations of $G$ up to isomorphism, the table of marks of $G$ classifies permutation 
representations of $G$ up to equivalence. Moreover, the table of marks encodes a wealth of information about the subgroup structure of $G$ in a compact way. For instance, 
up to a known factor, the mark $\beta_{G/H_{i}}(H_{j})$ is exactly the number of conjugates of the subgroup $H_{i}$ which contain $H_{j}$ as a subgroup. 

Thus, the table of marks provides a close approximation of the subgroup lattice of $G$ and precisely describes the poset of conjugacy classes of subgroups of $G$.
Conversely, the table of marks can be obtained by counting incidences in the subgroup lattice of $G$. However, both the computation of the subgroup lattice of $G$ as well as 
incidence counting between conjugacy classes of subgroups are computationally expensive tasks, unless the order of $G$ is small. It is therefore desirable to be 
able to compute the table of marks in a way that avoids computing the subgroup lattice, or counting incidences, or both.

Pfeiffer \cite{pfe} describes a procedure for the construction of the table of marks of a finite group $G$ from the tables of marks of its maximal subgroups. 
This semi-automatic procedure has proven well suited for simple groups up to a certain order, and has been used extensively in building the \gap\ \cite{GAP4}  library of tables 
of marks \tomlib~\cite{tomlib}.

In this article we present a new algorithm for the computation of the table of marks of a cyclic extension of $G$ from the table of marks of $G$. More 
precisely, we show how to compute the \emph{subgroup pattern} of the extension from the subgroup pattern of $G$. Here, the subgroup pattern (c.f.~\cite{geometriesandgroups,structurem12}) of a finite group $G$ 
is a list of representatives of its conjugacy classes of subgroups together with its table of marks.  As a motivating example we choose the symmetric 
group $S_{n}$ which contains the alternating group $A_{n}$ as a normal subgroup of index $2$. With this in mind, we will assume from Section~\ref{sec:subgroups} on that $S$ is a finite group, 
that $A$ is a normal subgroup of $S$ of index $p$ for some prime number $p$, and that the subgroup pattern of $A$ is known.    

In Section~\ref{sec:basics}, we introduce notation and review some
basic properties of $G$-sets and $G$-maps.  In
Section~\ref{sec:subgroups}, we describe an algorithm for the
computation of the conjugacy classes of subgroups of $S$ from a list
of representatives of the conjugacy classes of subgroups of~$A$.
Repeated application of this algorithm yields an algorithm for the
computation of the conjugacy classes of subgroups of a solvable group.
In Section~\ref{sec:marks}, we discuss the building blocks for the
computation of the table of marks of $S$ from the table of marks of
$A$, assuming that the conjugacy classes of subgroups of both $A$ and
$S$ are known.  In the final section, we combine these tools into
an algorithm for the computation of the subgroup pattern of $S$
from the subgroup pattern of $A$.
Repeated application of this algorithm yields an algorithm for the
computation of the table of marks of a solvable group.  The section
finishes with a list of concrete results and  performance statistics.



\section{$G$-sets and $G$-maps}\label{sec:basics}
Let $G$ be a finite group. A finite set $X$ together with a map $X \times G \to X$, mapping the pair $(x, g) \in X \times G $ to $x.g \in X$ is called a $G$-set 
if $x.1 = x \mbox{ for all } x \in X$ and $(x.g).g' = x.(gg') \mbox{ for all } x \in X,  g,g' \in G$. A map $f : X \to Y$ between $G$-sets $X$ and $Y$ is called a $G$-map 
if $f(x.g) = f(x).g \mbox{ for all } x \in X, g \in G$. We review some notation and basic properties of $G$-sets and the maps between them.

For a $G$-set $X$, we denote by $\pi_{X}:G \rightarrow \mathbb{N}_{0}$ the permutation character (see \cite{breuer}) of the action of $G$ on $X$, i.e.
\[ \pi_{X}(g) = |\mathrm{Fix}_{X}(g)| = \# \{ x \in X : x.g = g \}, \]
for $g \in G$. 

The group $G$ partitions any $G$-set $X$ into orbits. For $x \in X$, we denote by $[x]_{G} = x.G$ (or simply $[x]$) the $G$-orbit (or class) of $x$, and by 
\[ X/G = \{ [x]_{G} : x \in X \} \]
the quotient set (or set of classes). The number of orbits of $G$ on $X$ can be computed from the permutation character as
\begin{align}\label{burnlemma}
 |X/G| = \frac{1}{|G|} \sum_{g \in G} \pi_{X}(g),
\end{align}
by the Cauchy-Frobenius Lemma (the lemma that is not Burnside's \cite{neumann}).

If $G$ acts on two sets $X$ and $Y$ then $G$ also acts on their product $X \times Y$ via $(x, y).g = (x.g, y.g)$ for all $x \in X, y \in Y, g \in G$. 
The following propositions list some general properties of this action on pairs which will be used in the sequel. Their proofs make use of the following easy 
lemma.
\begin{lemma}
\label{pairs1}
 Suppose that $X$ and $Y$ are $G$-sets. Then,
\begin{enumerate}
 \item for all $x \in X$, $y \in Y$, we have
\[ [x, y]_{G} \cap (X \times {y}) = [x]_{G_{y}} \times {y}; \]
\item for $y \in Y$, the map $[x]_{G_{y}} \mapsto [x, y]_{G}$ is a well defined bijection from $X/G_{y}$ to $(X \times [y]_{G})/G$.
\end{enumerate}
\end{lemma}
\begin{proof}
(i) The statement is equivalent to 
\[ \{ x' \in X : (x', y) \in (x, y).G \} = x.G_{y} \]
which is obviously true.

(ii) Consider the map $\gamma : X \rightarrow (X \times Y)/G$ defined by $\gamma (x) = [x, y]_{G}$ for $x \in X$. Then 
$\gamma (X) = (X \times [y]_{G})/G$, and by (i), $\gamma^{-1}([x, y]_{G}) = [x]_{G_{y}}$.
\end{proof}
\begin{prop}\label{pairs2}
 Suppose that $X$ and $Y$ are transitive $G$-sets and that $Z \subseteq X \times Y$ is a $G$-invariant subset of pairs. Let $(x, y) \in Z$. Then the stabilizers 
$G_{y}, G_{x}$ act on 
\[ Zy = \{ x' \in X : (x', y) \in Z \}, \quad xZ = \{ y' \in Y : (x, y') \in Z \} \]
respectively, and the map $\xi : Zy/G_{y} \rightarrow xZ/G_{x}$, given by
\[ \xi ([x.a]_{G_{y}}) = [y.a^{-1}]_{G_{x}} \]
for $a \in G$, is a well defined bijection of orbits.
\end{prop}
\begin{proof}
 By Lemma \ref{pairs1}, the maps $\alpha : Zy/G_{y} \rightarrow Z/G$ and $\beta : xZ/G_{x} \rightarrow Z/G$, defined by 
\[ \alpha([x']_{G_{y}}) = [x',y]_{G}, \quad \beta ([y']_{G_{x}}) = [x,y']_{G} \]
\noindent for $x' \in Zy, y' \in xZ$, are well defined bijections, and $\xi = \beta^{-1} \circ \alpha$. 
\end{proof}
\begin{prop}\label{pairs3}
 Suppose that $X$ and $Y$ are $G$-sets and that $f : X \rightarrow Y$ is a $G$-map. Then the map
\[ \zeta : \coprod_{[y]\in Y/G} f^{-1}(y)/G_{y} \rightarrow X/G \]
defined by $\zeta ([x]_{G_{f(x)}}) = [x]_{G}$ for $x \in f^{-1}(y)$, where $y$ ranges over a set of representatives of the $G$-orbits on $Y$, is a 
well defined bijection.
\end{prop}
\begin{proof}
 The set $Z = \{ (x,y) \in X \times Y : y = f(x) \}$ is a $G$-invariant subset of $X \times Y$ with $xZ = \{ f(x) \}$ for all $x \in X$, and $Zy = f^{-1}(y)$ for all $y \in Y$, 
in the notation of Proposition \ref{pairs2}. By Lemma \ref{pairs1}, for each orbit $[y] \in Y/G$, there is a bijection $[x]_{G} \mapsto [x,y]_{G}$ between 
$f^{-1}(y)/G_{y}$ and $(X \times [y])/G$, which in turn is a bijection to $f^{-1}([y])/G$ via $[x, f(x)]_{G} \mapsto [x]_{G}$. The claim then follows from the fact 
that
\[ X = \coprod_{y \in Y} f^{-1}(y) = \coprod_{[y] \in Y/G} f^{-1}([y]_{G}) ,\]
whence $X/G = \coprod_{[y] \in Y/G} f^{-1}([y])/G$.
\end{proof}
\subsection{Marks}
We call the collection of all marks which $G$ leaves on $X$, that is the function $\beta_{X} : \mathrm{Sub}(G) \rightarrow \mathbb{Z} $, which assigns to each 
subgroup $H$ of $G$ its mark 
\[ \beta_{X}(H) = |\mathrm{Fix}_{X}(H)| = \# \{ x \in X : x.h = x \mbox{ for all } h \in H \}, \]
the \emph{impression} of $G$ on $X$. Clearly, $\beta_{X}$ is constant on conjugacy classes, so we can regard $\beta_{X}$ as a function from the set 
$\mathrm{Sub}(G)/G$ of conjugacy classes of subgroups of $G$ to $\mathbb{Z}$, or simply as the list of integers
\[ \beta_{X} = ( \beta_{X}(H_{1}), \ldots , \beta_{X}(H_{r}) )\] 
where $H_{1}, \ldots , H_{r}$ is a fixed list of representatives of the conjugacy classes of subgroups of $G$. The table of marks of $G$ is then the 
$r \times r$-matrix which has as its rows the impressions of the transitive $G$-sets $G/H_{i}, i = 1, \ldots, r$. Marks can also be viewed as incidences between 
conjugacy classes of subgroups due to the following formula (e.g., see~\cite[Prop 1.2]{pfe}):
\begin{align}\label{marksasincidence}
 \beta_{G/K}(H) = |N_{G}(K):K| \cdot \# \{K^g : H \leq K^g, g \in G \}.
\end{align}
\begin{theorem}[Burnside~\cite{burn}]
\label{burnsidestheorem}
Let $G$ be a finite group, and $X$ and $Y$ be finite $G$ sets. Then the $G$-sets $X$ and $Y$ are isomorphic if and only if $\beta_{X} = \beta_{Y}$.
\end{theorem}
\subsection{The Burnside Ring}
\label{burnsidering}
For any $G$-set $X$, let $[X]$ denote its isomorphism class. The Burnside ring of $G$, denoted $\Omega (G)$ is the free abelian group 
\[ \Omega (G) = \{ \sum \limits_{i = 1}^{r} a_{i} [G/H_{i}] : a_{i} \in \mathbb{Z} \} \] 
generated by the isomorphism classes of transitive $G$-sets $[G/H_{i}], i = 1, \ldots, r $. The sum $[X] + [Y]$ of the isomorphism classes of 
$G$-sets $X$ and $Y$ is the isomorphism class $ [X \sqcup Y]$ of the disjoint union of $X$ and $Y$, and the product $[X] \cdot [Y]$ is the isomorphism 
class  $[X \times Y]$ of the Cartesian product of $X$ and $Y$. This turns $\Omega (G)$ into a commutative ring with identity $[G/G]$ (see~\cite{bouc}). 
\subsection{Dress Congruences}\label{dressintro}
Note that, if $X$ and $Y$ are $G$-sets, and $H$ is a subgroup of $G$, then $\beta_{X \sqcup Y}(H) = \beta_{X}(H) + \beta_{Y}(H)$ and 
$\beta_{X \times Y}(H) = \beta_{X}(H) \times \beta_{Y}(H)$. Theorem \ref{burnsidestheorem} has the following consequence. Each subgroup 
$H$ of $G$ defines a ring homomorphism $\Omega (G) \rightarrow \mathbb{Z}$ by $[X] \mapsto \beta_{X}(H)$. Since 
$\beta_{X}(H) = \beta_{X}(K)$ if $H$ and $K$ are conjugate in $G$, it follows that the product mapping
\begin{eqnarray*}
 \beta : &\Omega (G)& \rightarrow \mathbb{Z}^{r} \\
 &[X]& \mapsto \beta_{X} = (\beta_{X}(H_{1}), \ldots , \beta_{X}(H_{r}) )
\end{eqnarray*}
is injective. In this context $\mathbb{Z}^{r}$ is often called the \emph{ghost ring} of $G$. 

The matrix $\mathrm{M}(G)$ of the linear map $\beta$ with respect to the basis $\{G/H_{i} \}_{i = 1, \ldots, r}$ of $\Omega(G)$ and to the canonical basis $\{ u_{i} \}_{i = 1, \ldots, r}$ of $\mathbb{Z}^r$ 
is the table of marks of $G$. Thus, if 
\[ [X] = \sum \limits_{i = 1}^{r} a_{i} [G/H_{i}] \in \Omega (G), \] 
then $\beta_{X}$ can be expressed 
in terms of the table of marks $\mathrm{M}(G)$ as 
\[ \beta_{X} = (a_{1}, \ldots, a_{r})  \mathrm{M}(G). \]
\begin{theorem}(Dress, see \cite{bouc, dress})
\label{dresscongtheorem}
  Let $G$ be a finite group. For $H, U \leq G$, set
\[ n(U,H) = \# \{ Ua \in N_{G}(U)/U : \langle U, a \rangle \sim_{G} H \}. \]
Then the element $y =  (y_{1},\ldots, y_{r})$ of $\mathbb{Z}^{r}$ is in the image of $\beta$ if and 
only if 
\[ \sum \limits_{i=1}^{r} n(U,H_{i}) y_{i} \equiv 0 \mbox{ mod } |N_{G}(U)/U|. \]
for all $U \leq G$.
\end{theorem}
Theorem \ref{dresscongtheorem} yields a set of congruences which, in particular, must be satisfied by the rows of the table of marks of $G$.

\section{The Subgroups of $S$}\label{sec:subgroups}
\begin{figure}[H]

  \begin{center}
\scalebox{0.8}{
\begin{tikzpicture}

 \node (1) at (10,0)[blue]   {$1$};
\node (b2) at (8,2)[blue]   {$2$};

\node (r2) at (6,2)[red]  {$2$};
\node (b3) at (12, 3)[blue]   {$3$};
\node (b22) at (8, 4)[blue]   {$2^2$};
\node (r22) at (4, 4)[red]  {$2^2$};
\node (r4) at (6, 4)[red]  {$4$};
\node (s3) at (13, 5)[red]  {$S_{3}$};
\node (d8) at (6, 6) [red] {$D_{8}$};
\node (a4) at (10, 8) [blue] {$A_{4}$};
\node (s4) at (10, 10)[red]  {$S_{4}$};

\path
(1) edge[thick, blue] (b2)
(1) edge[thick,blue]  (b3)
(b2) edge[thick,blue]  (b22)
(b3) edge[thick,blue]  (a4)
(b22) edge[thick,blue]  (a4)
(1) edge[thick, red, dashed] (r2)
(b2) edge[thick,red, dashed] (r22)
(b2) edge[thick,red, dashed] (r4)
(b22) edge[thick,red, dashed] (d8)
(a4) edge[thick, red, dashed] (s4)
(b3) edge[thick,red, dashed] (s3)
(r2) edge[thick] (s3)
(r2) edge[thick]  (r22)
(r4) edge[thick]  (d8)
(r22) edge[thick]  (d8)
(d8) edge[thick]  (s4)
(s3) edge[thick]  (s4);

\end{tikzpicture}}
\end{center}
\caption{Poset of Conjugacy Classes of Subgroups of $S_{4}$}
\label{lattices4}
 \end{figure}
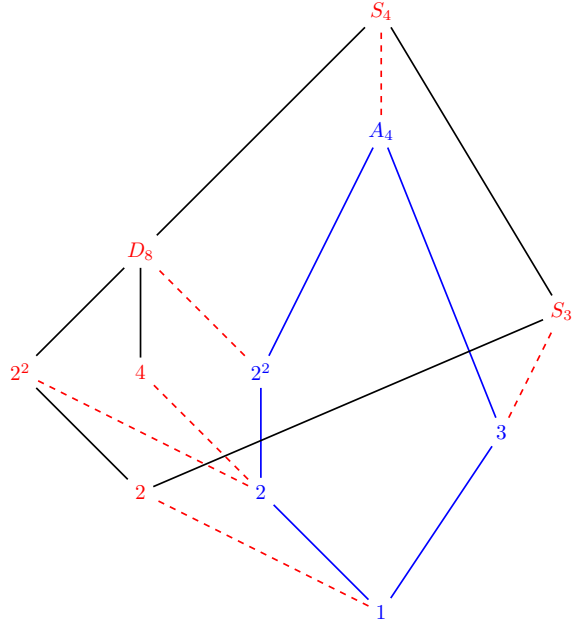

From now on, let $S$ be a finite group, and let $A$ be a normal subgroup of $S$ of index $p$ for some prime $p$.
In this section we describe an algorithm for the computation of the conjugacy classes of subgroups of $S$ from the conjugacy classes of 
subgroups of~$A$. For the purpose of exposition we distinguish between two types of subgroups of~$S$: the subgroups of $A$ will be called \textit{blue subgroups}, and the subgroups of $S$ which are not contained 
in $A$ will be called \textit{red subgroups}. The set of subgroups of $S$ then is a disjoint union 
\[ \mathrm{Sub}(S) = \mathcal{B} \sqcup \mathcal{R}, \]
where 
\[ \mathcal{B} = \mathrm{Sub}(A),\quad \mathcal{R} = \mathrm{Sub}(S) \setminus \mathrm{Sub}(A). \]
Since no red subgroup is conjugate to a blue subgroup, both $\mathcal{B}$ and $\mathcal{R}$ are $S$-sets. The aim of this section is 
to obtain an effective description of the conjugacy classes  \[ \mathrm{Sub}(S)/S = \mathcal{B}/S \sqcup \mathcal{R}/S \] of subgroups of $S$ from the 
conjugacy classes $\mathrm{Sub}(A)/A = \mathcal{B}/A$ of subgroups of $A$. As a simple example, the separation of $\mathrm{Sub}(S_{4})/S_{4} $ into blue and red classes of subgroups is illustrated 
in Figure \ref{lattices4}, where, blue subgroups are connected by blue edges, red subgroups are connected by black edges, and dashed red eges are used to connect 
blue subgroups to red subgroups.

\subsection{Classes of Blue Subgroups}
Blue conjugacy classes of subgroups  of $S$ are unions of $A$-conjugacy classes of subgroups of $A$. The following proposition shows that a blue conjugacy class 
in $S$ is a union of exactly one or $p$ $A$-conjugacy classes.
\begin{prop}\label{bluesubs}
 Let $H \leq A$ and let $t \in S\setminus A$. Then 
\[[H]_{S} = [H]_{A} \cup [H^t]_{A} \cup \cdots \cup [H^{t^{p^{-1}}}]_{A} \]
where either $[H]_{S} = [H]_{A}$ and $| N_{S}(H) : N_{A}(H) | = p$, or $N_{S}(H) = N_{A}(H)$ and $|[H]_{S}| = p \left| [H]_{A}\right|$.
\end{prop}
\begin{proof}
 First, note that each $S$-conjugate of $H$ lies in one of $[H]_{A}, [H^t]_{A}, \ldots , [H^{t^{p^{-1}}}]_{A}$, since 
$S = A \cup tA\cup \cdots \cup t^{p^{-1}}A$. Moreover, each of the $A$-conjugacy classes of the $S$-conjugates of $H$ have the same size, since conjugation by $t$ 
induces a bijection between $[H]_{A}$ and $[H^t]_{A}$. By the Orbit-Stabilizer Theorem,
\[ |[H]_{S}|\cdot|N_{S}(H)| = |S| = p \left|A \right| = p \left|[H]_{A}\right|\cdot|N_{A}(H)| \]
From $[H]_{A} \subseteq [H]_{S}$ and $N_{A}(H) \leq N_{S}(H)$, it follows that either $[H]_{A} = [H]_{S}$ and $|N_{S}(H)| = p \left|N_{A}(H)\right|$ or that $N_{S}(H) = N_{A}(H)$ 
and $|[H]_{S}| = p \left|[H]_{A}\right|$.
\end{proof}
According to the dichotomy in this proposition, we denote
\[ \mathcal{B}_{1} = \{ H \in \mathcal{B}: [H]_{S} = [H]_{A} \}, \mbox{    } \mathcal{B}_{2} = \{ H \in \mathcal{B} : N_{S}(H) = N_{A}(H) \} \]
Then $\mathcal{B} =  \mathcal{B}_{1} \sqcup \mathcal{B}_{2}$ implies $\mathcal{B}/A =  \mathcal{B}_{1}/A \sqcup \mathcal{B}_{2}/A$ and 
the $S$-conjugacy class of blue subgroups can be described as follows.
\begin{cor}\label{numsubs}
 $\mathcal{B}/S =  \mathcal{B}_{1}/A  \sqcup  \mathcal{B}_{2}/S$. In particular, $S$ has $b = b_{1} + \frac{1}{p}b_{2}$ conjugacy classes of blue subgroups, 
where $b_{i} = |\mathcal{B}_{i}/A|, i = 1,2$.
\end{cor}
Corollary \ref{numsubs} yields the following algorithm to compute the set $\mathcal{B}/S$ of blue subgroups of $S$ from the set $\mathcal{B}/A$.

\begin{algorithm}[H]   
\caption{ \texttt{BlueSubgroups()} }
\label{blue}                           
\begin{algorithmic}
\State \textbf{Input} Representatives of $\mathcal{B}/A$
\State \textbf{Output} Representatives of $\mathcal{B}/S$
\State Initialize $B_{1} \gets\{ \},\, B_{2} \gets \{ \}$
\For{ $H \in \mathcal{B}/A$} 
\If{ $N_{S} (H) \nleq A$}
\State Add $H$ to $B_{1}$.
\Else
\State Add $H$ to $B_{2}$.
\EndIf
\EndFor
\State \textbf{return} $B_{1} \cup$ (a set of representatives of $S$-conjugate subgroups in $B_{2}$).

\end{algorithmic}
\end{algorithm}

\begin{example}
The special linear group $L_{2}(32)$ is a normal subgroup of index $5$ in $L_{2}(32){:}5$. Figure \ref{l232} illustrates how the blue classes of subgroups of $L_{2}(32)$ 
fuse to form blue classes of subgroups of $L_{2}(32){:}5$.
\begin{figure}[H]
\scalebox{0.8}{

\begin{tikzpicture}

\node (l232) at (0,0) {\small{$L_{2}(32)$}};
\node (l1) at (1,0) {\small{$1$}};
\node (l2) at (1.5, 0) {\small{$C_{2}$}};
\node (l3) at (2,0) {\small{$C_{3}$}};
\node (l221) at (2.5, 0) {\small{$2^2$}};
\node (l222) at (2.9, 0) {\small{$2^2$}};
\node (l223) at (3.3, 0) {\small{$2^2$}};
\node (l224) at (3.7, 0) {\small{$2^2$}};
\node (l225) at (4.1, 0) {\small{$2^2$}};
\node (ls3) at (4.5,0) {\small{$S_{3}$}};
\node(l2221) at (5,0) {\small{$2^3$}};
\node(l2222) at (5.4,0) {\small{$2^3$}};
\node(l2223) at (5.8,0) {\small{$2^3$}};
\node(l2224) at (6.2,0) {\small{$2^3$}};
\node(l2225) at (6.6,0) {\small{$2^3$}};
\node (lc11) at (7.1, 0) {\small{$C_{11}$}};
\node (l2^4) at (7.6,0) {\small{$2^4$}};
\node (ld22) at (8.1, 0) {\small{$D_{22}$}};
\node (lc31) at (8.7,0) {\small{$C_{31}$}};
\node (l2^5) at (9.2,0) {\small{$2^5$}};
\node (lc33) at (9.7,0) {\small{$C_{33}$}};
\node (d62) at (10.3,0) {\small{$D_{62}$}};
\node (ld66) at (10.9,0) {\small{$D_{66}$}};
\node (split) at (11.9,0) {\small{$2^5{:}C_{31}$}};
\node (l) at (13,0) {\small{$L_{2}(32)$}};

\node (5l232) at (0,4) {\small{$L_{2}(32){:}5$}};
\node (5l1) at (1,4) {\small{$1$}};
\node (5l2) at (1.5, 4) {\small{$C_{2}$}};
\node (5l3) at (2,4) {\small{$C_{3}$}};
\node (5l221) at (3.3, 4) {\small{$2^2$}};
\node (5ls3) at (4.5,4) {\small{$S_{3}$}};
\node(5l2221) at (5.8,4) {\small{$2^3$}};
\node (5lc11) at (7.1, 4) {\small{$C_{11}$}};
\node (5l2^4) at (7.6,4) {\small{$2^4$}};
\node (5ld22) at (8.1, 4) {\small{$D_{22}$}};
\node (5lc31) at (8.7,4) {\small{$C_{31}$}};
\node (5l2^5) at (9.2,4) {\small{$2^5$}};
\node (5lc33) at (9.7,4) {\small{$C_{33}$}};
\node (5d62) at (10.3,4) {\small{$D_{62}$}};
\node (5ld66) at (10.9,4) {\small{$D_{66}$}};
\node (5split) at (11.9,4) {\small{$2^5{:}C_{31}$}};
\node (5l) at (13,4) {\small{$L_{2}(32)$}};

\path

(l1) edge[thick, blue] (5l1)
(l2) edge[thick, blue] (5l2)
(l3) edge[thick, blue] (5l3)
(l221) edge[thick, blue] (5l221)
(l222) edge[thick, blue] (5l221)
(l223) edge[thick, blue] (5l221)
(l224) edge[thick, blue] (5l221)
(l225) edge[thick, blue] (5l221)
(ls3) edge[thick, blue] (5ls3)
(l2221) edge[thick, blue] (5l2221)
(l2222) edge[thick, blue] (5l2221)
(l2223) edge[thick, blue] (5l2221)
(l2224) edge[thick, blue] (5l2221)
(l2225) edge[thick, blue] (5l2221)
(lc11) edge[thick, blue] (5lc11)
(l2^4) edge[thick, blue] (5l2^4)
(ld22) edge[thick, blue] (5ld22)
(lc31) edge[thick, blue] (5lc31)
(l2^5) edge[thick, blue] (5l2^5)
(lc33) edge[thick, blue] (5lc33)
(d62) edge[thick, blue](5d62)
(ld66) edge[thick, blue] (5ld66)
(split) edge[thick, blue] (5split)
(l) edge[thick, blue] (5l);

\end{tikzpicture}}
\caption{Class Fusions in $L_{2}(32){:}5$}
\label{l232}
\end{figure}
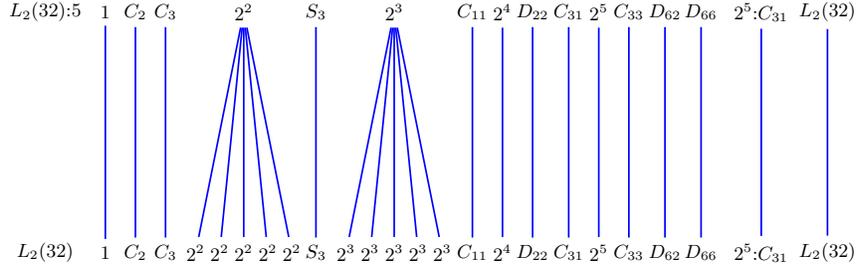
\end{example}
\subsection{Classes of Red Subgroups}
Red conjugacy classes of subgroups of $S$ correspond to certain conjugacy classes of subgroups of order $p$ in normalizer quotients.
\begin{prop}\label{redsubs}
 For $H \in \mathcal{B}$, let $T_{H} \subseteq S$ be such that $\{ H\langle t \rangle : t \in T_{H} \}$ is a transversal of the conjugacy classes of 
subgroups of order $p$ of $N_{S}(H)/H$ which lie outside $N_{A}(H)/H$. Then the set
\[ \coprod \limits_{[H]_{A} \in \mathcal{B}/A}^{} \{ \langle H, t \rangle : t \in T_{H} \}, \]
where $H$ ranges over a transversal of $\mathcal{B}/A$, is a transversal of $\mathcal{R}/S$.
\end{prop}
\begin{proof}
 Consider the map $\gamma : \mathcal{R} \rightarrow \mathcal{B}$, defined by $\gamma (K) = A \cap K$ for $K \in \mathcal{R}$. From
\[ \gamma (K^{s}) = K^s \cap A = K^s \cap A^s = (K \cap A)^s = \gamma(K)^s \]
for any $s \in S$, it follows that $\gamma$ is an $S$-map. For $H \in \mathcal{B}$, the map $K \mapsto K/H$ is a bijection between
\[ \mathcal{R}_{H} = \{ K \in \mathcal{R} : \gamma (K) = H \} = \gamma^{-1}(H) \]
and the set of subgroups of order $p$ in the quotient $N_{S}(H)/H$ which are not contained in $N_{A}(H)/H$. Moreover, these two sets are equivalent as 
$N_{S}(H)/H$-sets. By Proposition \ref{pairs3}, 
\[ \mathcal{R}/S = \coprod \limits_{[H]_{S} \in \mathcal{B}/S}^{} \mathcal{R}_{H}/N_{S}(H), \]
where $H$ ranges over a transversal of the conjugacy classes of blue subgroups of $S$. The statement remains true, if $H$ ranges over a transversal of 
$\mathcal{B}_{1}/S = \mathcal{B}_{1}/A$, or over a transversal of $\mathcal{B}/A$, since $\mathcal{R}_{H} = \emptyset$ for all $H \in \mathcal{B}_{2}$.
\end{proof}

Note that $T_{H} \subseteq S$ can easily be determined from a list of representatives of the conjugacy classes of $N_{S}(H)/H$. In fact, modulo $H$, the set $T_{H}$ is in 
bijection to the set of rational classes of elements of order $p$ in $N_{S}(H)/H \setminus N_{A}(H)/H$.
Moreover, each $t \in T_{H}$ can be chosen to 
be an element of order a power of $p$. 
\begin{cor}\label{numred}
 With the above notation, $S$ has 
\[ r = \sum \limits_{[H]_{A} \in \mathcal{B}/A} |\mathcal{R}_{H}/N_{S}(H)| = \sum \limits_{[H]_{A} \in \mathcal{B}/A} |T_{H}| \] 
conjugacy classes of red subgroups.
\end{cor}
Proposition \ref{redsubs} yields the following algorithm to compute the set $\mathcal{R}/S$ of red subgroups of $S$.

\begin{algorithm}[H]                 
\caption{ \texttt{RedSubgroups()}}  
\label{red}        
\begin{algorithmic}[H]
\State \textbf{Input} Representatives of $\mathcal{B}/A$
\State \textbf{Output} Representatives of $\mathcal{R}/S$
\State output $\gets \{ \}$.
\For{$H \in \mathcal{B}/A$}
\If {$N_{S} (H) \nleq A$}
\State Use \texttt{RationalClasses(}$N_{S}(H)/H$\texttt{)} to compute $T_{H}$
\For{$t \in T_{H}$}
\State Append $\{ \langle H, t \rangle : t \in T_{H} \}$ to output. 
\EndFor
\EndIf
\EndFor
\State \textbf{return} $\mathcal{R}/S$.
\end{algorithmic}
\end{algorithm}

\noindent It follows with Corollaries \ref{numsubs} and \ref{numred} that $|\mathrm{Sub}(S)/S| = b + r$. The $b + r$ conjugacy classes of subgroups of $S$ can 
now be enumerated by the following combination of Algorithms \ref{blue} and \ref{red}.

\begin{algorithm}[H]   
\caption{\texttt{SubgroupsByCyclicExtension()}} \label{subextensions}
\begin{algorithmic}[H]
\State \textbf{Input} Representatives of $\mathcal{B}/A$.
\State \textbf{Output} Representatives of $\mathrm{Sub}(S)/S$.
\State \textbf{return} \texttt{BlueSubgroups(}$\mathcal{B}/A$\texttt{)} $\cup$ \texttt{RedSubgroups(}$\mathcal{B}/A$\texttt{)}.
\end{algorithmic}
\end{algorithm}
Recall from the introduction that the subgroup pattern of $S$ consists of the list of representatives of the conjugacy classes of 
subgroups of $S$ and the table of marks of $S$. Accordingly, the task of computing the subgroup pattern of $S$ from that of $A$
requires the computation of the conjugacy classes of subgroups of $S$ from those of $A$, and the computation of the table of marks 
of $S$ from that of $A$. Algorithm \ref{subextensions} accomplishes the first part of this task.

\subsection{Computing the Subgroups of a Solvable Group}\label{solvablesubs}
Algorithm \ref{subextensions} has enabled us to produce a new algorithm to compute the conjugacy classes of subgroups of a solvable 
group $G$ in an iterative fashion starting with the conjugacy classes of subgroups of the trivial group. Recall that a solvable group $G$  has a
composition series of the form 
\[ 1 = G_{0} \trianglelefteq G_{1} \ldots \trianglelefteq G_{n} = G \]
in which each factor $G_{i + 1} / G_{i}$ is cyclic of prime order. In such cases we can apply the methods described in Propositions \ref{bluesubs} and 
\ref{redsubs} to compute the conjugacy classes of subgroups of $G$ in a step by step fashion.

\begin{algorithm}[H]                 
\caption{ \texttt{AllSubgroupClassesSolvable()}}          
\begin{algorithmic}[H]
\State \textbf{Input} A solvable group $G$.
\State \textbf{Output} $\mathrm{Sub}(G)/G$.
\State Compute a composition series $1 = G_{0} \trianglelefteq G_{1} \trianglelefteq \ldots \trianglelefteq G_{n} = G$
\State Obviously $\mathrm{Sub}(G_{0}) = \{ 1 \}$.
\For{ $i \in \{1, \ldots, n \}$ }
\State Compute $\mathrm{Sub}(G_{i})/G_{i}$ as \texttt{SubgroupsByCyclicExtension(}$\mathrm{Sub}(G_{i-1})/G_{i-1}$\texttt{)}.
\EndFor
\State \textbf{return} $\mathrm{Sub}(G)/G$.
\end{algorithmic}
\end{algorithm}

The performance of our implementation of this algorithm in \texttt{GAP} compares quite favourably to the existing \texttt{GAP} functions for computing conjugacy classes of subgroups, notably 
\texttt{SubgroupsSolvableGroup} (see \cite{hulpke}), and the standard \texttt{GAP} function \texttt{ConjugacyClassesSubgroups} for computing conjugacy classes of subgroups.

\begin{example}
 Consider the General linear group GL$_{2} (3)$ of all invertible $2 \times 2$ matrices over the field with $3$ elements. $\mathrm{GL}_{2} (3)$ is a solvable 
group and has the following composition series
\[1  \vartriangleleft 2 \vartriangleleft 4 \vartriangleleft Q_{8} \vartriangleleft \textrm{SL}_{2} (3) \vartriangleleft \mathrm{GL}_{2} (3) \]

\noindent Figure \ref{gl23} shows the growth and fusion of conjugacy classes of subgroups as we incrementally extend from one group in the 
composition series to the next. 

\begin{figure}[H]
\scalebox{0.8}{
\begin{tikzpicture}

\node (1) at (0,0) {$1$};
\node (2) at (0,2) {$2$};
\node (4) at (0,4) {$4$};
\node (q8) at (0,6) {$Q_{8}$};
\node (sl) at  (0,8) {$\mathrm{SL}_{2}(3)$};
\node (gl) at (0,10) {$\mathrm{GL}_{2}(3)$};

\path (1) edge[thick] (2)
(2) edge[thick] (4)
(4) edge[thick] (q8)
(q8) edge[thick] (sl)
(sl) edge[thick] (gl);

\node (11) at (2, 0) {$1$};
\node (21) at (2, 2) {$1$};
\node (22) at (3,2) {$2$};
\path (11) edge[thick, blue] (21)
      (11) edge[thick, dashed,red] (22);

\node (41) at (2,4) {$1$};
\node (42) at (3,4) {$2$};
\node (44) at (3.5,4) {$4$};
\path (41) edge[thick, blue] (21)
      (22) edge[thick, blue] (42)
      (22) edge[thick, dashed,red] (44);

\node (q1) at (2,6) {$1$};
\node (q2) at (3,6) {$2$};
\node (q41) at (3.5,6) {$4$};
\node (q42) at (4.5,6)  {$4$};
\node (q43) at (5.5,6) {$4$};
\node (q8) at  (8,6) {$Q_{8}$};

\path 
(41) edge[thick, blue] (q1)
(42) edge[thick, blue] (q2)
(44) edge[thick, blue] (q43)
(42) edge[thick, dashed,red] (q41)
(42) edge[thick, dashed,red] (q42)
(44) edge[thick, dashed,red] (q8);

\node (s1) at (2,8) {$1$};
\node (s2) at (3,8) {$2$};
\node (s3) at (3.5,8) {$3$};
\node (s4) at (4.5,8) {$4$}; 
\node (s6) at (6,8) {$6$};
\node (s8) at (8, 8) {$Q_{8}$};
\node (sl23) at (11, 8) {$\mathrm{SL}_{2}(3)$};
\path
(q1) edge[thick, dashed,red] (s3)
(q2) edge[thick, dashed,red] (s6)
(q8) edge[thick, dashed,red] (sl23)
(q1) edge[thick, blue] (s1)
(q2) edge[thick, blue] (s2)
(q41) edge[thick, blue] (s4)
(q42) edge[thick, blue] (s4)
(q43) edge[thick, blue] (s4)
(q8) edge[thick, blue] (s8);

\node (g1) at (2,10) {$1$};
\node (g21) at (3,10) {$2$};
\node (g22) at (2.5,10) {$2$};
\node (g3) at (3.5,10) {$3$};
\node (g4) at (4.5,10) {$4$};
\node (g2sq) at (4,10) {$2^2$};
\node (gs31) at (5,10) {$S_{3}$};
\node (gs32) at (5.5,10) {$S_{3}$};
\node (g6) at (6,10) {$6$};
\node (gq8) at (8,10) {$Q_{8}$};
\node (gd8) at (6.8,10) {$D_{8}$};
\node (g8) at (7.4,10) {$8$};
\node (gd12) at (9.2,10) {$D_{12}$};
\node (gqd16) at (9.9,10) {$16$};
\node (gsl) at (11,10) {$\mathrm{SL}_{2}(3)$};
\node (glgl) at (12.5,10) {$\mathrm{GL}_{2}(3)$};

\path 
(s1) edge[thick, blue] (g1)
(s2) edge[thick, blue] (g21)
(s3) edge[thick, blue] (g3)
(s4) edge[thick, blue] (g4)
(s6) edge[thick, blue] (g6)
(s8) edge[thick, blue] (gq8)
(sl23) edge[thick, blue] (gsl)
(s1) edge[thick, dashed,red] (g22)
(s2) edge[thick, dashed,red] (g2sq)
(s3) edge[thick, dashed,red] (gs31)
(s3) edge[thick, dashed,red] (gs32)
(s4) edge[thick, dashed,red] (gd8)
(s4) edge[thick, dashed,red] (g8)
(s6) edge[thick, dashed,red] (gd12)
(s8) edge[thick, dashed,red] (gqd16)
(sl23) edge[thick, dashed,red] (glgl);

\end{tikzpicture}}
\caption{Class Fusions in $\mathrm{GL}_{2}(3)$}
\label{gl23}
\end{figure}
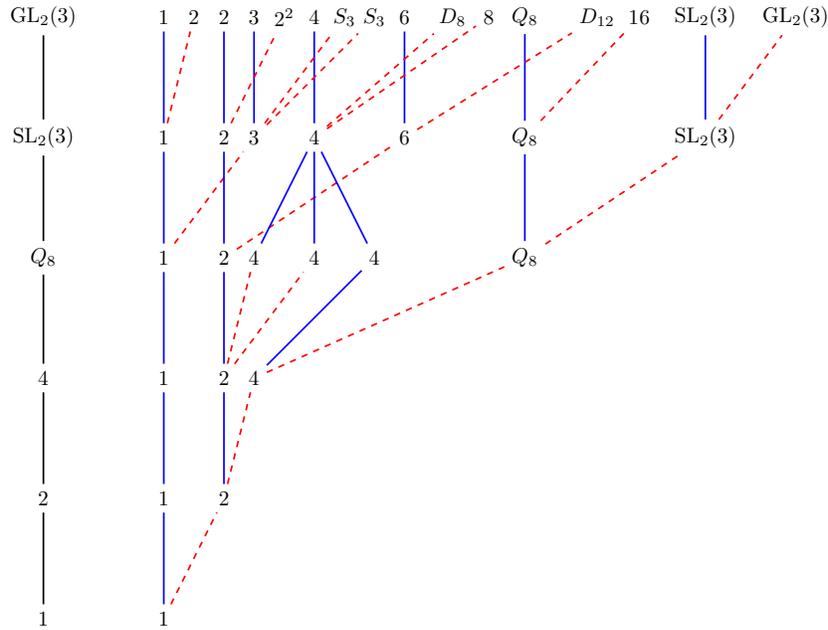
\end{example}

\section{The Table of Marks of $S$}\label{sec:marks}
In this section we develop tools for the computation of the table of marks of $S$ from the table of marks of $A$. For the purpose of describing the 
table of marks of $S$ in terms of the table of marks of $A$, we use the partition of the subgroups of $S$ into blue and red subgroups to subdivide the 
table of marks of $S$ into four quarters, labeled by pairs of colors. We illustrate the situation with the example of the alternating group $A_{5}$ as a 
subgroup of index $p = 2$ of the symmetric group $S_{5}$. The table of marks of $A_{5}$ is shown Figure \ref{toma5}.
  
\begin{center}
\begin{figure}[H]
\[\begin{array}{l|rrrrrrrrr}

A_{5}/1&60 & × & × & × & × & × & × & × & ×\\
A_{5}/C_{2}&30 & 2 & × & × & × & × & × & × & ×\\
A_{5}/C_{3}&20 & . & 2 & × & × & × & × & × & ×\\
A_{5}/2^2&15 & 3 & . & 3 & × & × & × & × & ×\\
A_{5}/C_{5}&12 & . & . & . & 2 & × & × & × & ×\\
A_{5}/S_{3}&10 & 2 & 1 & . & . & 1 & × & × & ×\\
A_{5}/D_{10}&6 & 2 & . & . & 1 & . & 1 & × & ×\\
A_{5}/A_{4}&5 & 1 & 2 & 1 & . & . & . & 1 & ×\\
A_{5}/A_{5}&1 & 1 & 1 & 1 & 1 & 1 & 1 & 1 & 1\\
\hline
&1&C_{2}&C_{3}&2^2&C_{5}&S_{3}&D_{10}&A_{4}&A_{5}\\
\end{array}\]
\caption{Table of Marks of $A_{5}$}
\label{toma5}
\end{figure}
\end{center}

The subdivided table of marks of $S_{5}$ is shown in Figure \ref{subdivideds5}.
\begin{center}

\begin{figure}[H]
\arraycolsep 2pt
$\begin{array}{l|rrrrrrrrr|rrrrrrrrrr}

S_{5}/1&120 & &  & × & & × & × & × & × &  &  &  &  &  & × & × & × & × & ×\\

S_{5}/C_{2}&60 & 4 & × &  & × & × & × & × & × &  &  &  & × & × & × & × & × & × & ×\\

S_{5}/C_{3}& 40 & . & 4 & × &  & & × & × & × & × & × &  &  & & × & × & × & × & ×\\

S_{5}/2^2&30 & 6 & . & 6 & × & × &  &  & × & × & × & × & × &  &  & × & × & × & ×\\

S_{5}/C_{5}&24 & . & . & . & 4 & × &  & × & × & × & × & × & × & × & × & × &  & × & ×\\

S_{5}/S_{3}&20 & 4 & 2 & . & . & 2 & × & × & × & × & × & × & × & × & × & & × & × & ×\\

S_{5}/D_{10}&12 & 4 & . & . & 2 & . & 2 & × & × & × & × & × & × & × & × & × &  & × & ×\\

S_{5}/A_{4} &10 & 2 & 4 & 2 & . & . & . & 2 & × & × & × & × & × & × & × & × & × &  & ×\\

S_{5}/A_{5}&2 & 2 & 2 & 2 & 2 & 2 & 2 & 2 & 2 & × & × & × & × & × & × & × & × & × & \\
\hline

S_{5}/C_{2} &60 & . & . & . & . & . & . & . & . & 6 & × &  & × &  & × & × & × & × & ×\\

S_{5}/C_{4}&30 & 2 & . & . & . & . & . & . & . & . & 2 & × & × & × &  & × & × & × & ×\\

S_{5}/2^2&30 & 2 & . & . & . & . & . & . & . & 6 & . & 2 & × & × &  & × & × & × & ×\\

S_{5}/S_{3}&20 & . & 2 & . & . & . & . & . & . & 6 & . & . & 2 & × & × & &  & × & ×\\

S_{5}/C_{6}&20 & . & 2 & . & . & . & . & . & . & 2 & . & . & . & 2 & × &  & × & × & ×\\

S_{5}/D_{8}&15 & 3 & . & 3 & . & . & . & . & . & 3 & 1 & 1 &  .& . &1  & × & × & × & ×\\

S_{5}/D_{12}&10 & 2 & 1 & . & . & 1 & . & . & . & 4 & . & 2 & 1 & 1 & . & 1 & × & × & ×\\

S_{5}/5{:}4&6 & 2 & . & . & 1 & . & 1 & . & . & . & 2 & . & . & . & . & . & 1 & × & ×\\

S_{5}/S_{4}&5 & 1 & 2 & 1 & . & . & . & 1 & . & 3 & 1 & 1 & 2 & . & 1 & . & . &  1& ×\\

S_{5}/S_{5}&1 & 1& 1& 1 &1 & 1& 1& 1 &1 & 1& 1& 1 &1 & 1& 1&1  &1 &1 & 1 \\
\hline
& 1& C_{2}& C_{3}& 2^2& C_{5}& S_{3}& D_{10}& A_{4}& A_{5}& C_{2}& C_{4}& 2^2& S_{3}& C_{6}& D_{8}& D_{12}& 5{:}4& S_{4}& S_{5} \\

\end{array}$

\caption{Table of Marks of $S_{5}$}
\label{subdivideds5}
\end{figure}
\end{center}
Since no red subgroup can be contained in any blue subgroup, the top right quarter, which represents the fixed points of red subgroups on blue subgroups, is zero. 
In this example, the top left quarter, which represents the fixed points of blue subgroups on blue subgroups, is exactly $p$ times the table of marks of 
$A_{5}$. The bottom left quarter, which represents the fixed points of blue subgroups on red subgroups, looks like a modified copy of the table of marks of $A_{5}$, in the 
sense that some rows are repeated, and the row in the table of marks of $A_{5}$ corresponding to $A_{5}/C_{5}$ does not appear at all. The bottom right quarter, 
which represents the fixed points of red subgroups on red subgroups, does not bear any immediate resemblance to the table of marks of $A_{5}$.
In the following sections we 
will examine each of these nonzero quarters separately.

\subsection{The Top Left Quarter}
Recall that the marks in this quarter represent fixed points of blue groups on blue groups. 
\begin{prop}
\label{blueonblue}
 Suppose that $H, U \leq A$. let $t \in S \setminus A$ and denote $H_{j} = H^{t^{j}}$, for $j = 0,1,\ldots ,p-1$. Then
\[ \beta_{S/H}(U) = \sum \limits_{j=0}^{p-1} \beta_{A/H_{j}}(U). \]
In particular if $[H]_{S} = [H]_{A}$ then $ \beta_{S/H} (U) = p  \beta_{A/H}(U)$.
 
\end{prop}
\begin{proof}
 The coset space $S/H$ is a disjoint union of $U$-sets, $\{ Hat^{j} : a \in A \} = \{Ht^{j}a : a \in A \}$ equivalent to 
$A/H_{j} = \{ H^{t^{j}}a : a \in A \}$, for $j = 0,1,\ldots , p-1$.
\end{proof}
If $\mathcal{B}_{2} = \emptyset$ then Proposition \ref{blueonblue} implies that the top left quarter of the table of marks of $S$ will be exactly $p$ times the 
table of marks of $A$ as observed in the example of $A_{5}$ and $S_{5}$. In general this quarter has one row for each class $[H]$ in $\mathcal{B}/S = \mathcal{B}_{1}/S \sqcup \mathcal{B}_{2}/S$, where, 
if $[H] \in \mathcal{B}_{1}/S$ the row is a $p$-multiple of the corresponding row in the table of marks of $A$, and if $[H] \in \mathcal{B}_{2}/S$ the row is then the sum 
of the rows corresponding to the $p$ $A$-conjugacy classes of subgroups which fuse to form a single $S$-conjugacy class of subgroups.
\subsection{The Bottom Left Quarter}\label{sec:bottomrightquarter}
Recall that the marks in this quarter represent the fixed points of blue subgroups on red subgroups.
\begin{prop}\label{blueonred}
 Suppose that $K \leq S$ is a red subgroup with $\gamma(K) = H \leq A$. then the coset spaces $S/K$ and $A/H$ are equivalent as $A$-sets. In particular,
\[ \beta_{S/K}(U) = \beta_{A/H}(U) \]
for all subgroups $U \leq A$.
\end{prop}
\begin{proof}
 The map $f : A/H \rightarrow S/K$, defined by $f(Ha) \mapsto Ka$ for $a \in A$, is an $A$-equivariant bijection and thus the coset spaces are equivalent as $U$-sets 
as well.
\end{proof}
It follows that for any $K \in \mathcal{R}$ with $\gamma (K) = H$ we insert a copy of the row corresponding to $H$ in the table of 
marks of $A$ into the bottom left quarter of the table of marks of $S$. This accounts for the duplicate rows observed in the example of $A_{5}$ and $S_{5}$.  
\subsection{The Bottom Right Quarter}
Recall that the marks in the bottom right quarter represent the fixed points of red subgroups on red subgroups. The marks in this section usually cannot be computed 
from the table of marks of $A$ using a simple formula. There are, however, obvious lower and upper bounds on these numbers, and various conditions which reduce the 
number of values that a particular mark can take. If a mark is not uniquely determined by these conditions, one can still compute it explicitly by counting incidences 
between the relevant conjugacy classes of subgroups. In this section we describe these bounds and conditions on the marks in question and describe how they can be completely determined. 
\subsubsection{Bounds}
\label{bounds}
 The marks in the bottom left quarter yield a first upper bound for the marks in the bottom right quarter.
\begin{lem}\label{ub}
 Let $H \leq K \leq S$. Then 
\[ \beta_{S/U} (K) \leq \beta_{S/U} (H) \]
for all subgroups $U \leq S$.
\end{lem}
\begin{proof}
Since $H \leq K$, clearly $K$ cannot fix more cosets than $H$.
\end{proof}

In particular if $K$ is a red subgroup with $\gamma (K) = H \leq A$ then $\beta_{S/U} (K) \leq \beta_{S/U}(H)$. Thus the marks in the bottom left quarter, provide 
an upper bound for the marks in the bottom right quarter. Combining Lemma \ref{ub} with the following Proposition we obtain a finite range of values for each of the 
marks in the bottom right quarter.
\begin{lem}
\label{congruentcolumns}
 Suppose $U, V \leq S$ with $U \trianglelefteq V$ of index $q$ a prime, and let $X$ be an $S$-set. Then 
\[ \beta_{X}(U) \equiv \beta_{X}(V) \mbox{ mod } q. \]
\end{lem}
\begin{proof}
 Clearly, $\mathrm{Fix}_{X}(U)$ can be regarded as a $V/U$-set. Since the quotient $V/U$ is cyclic of prime order, It follows that $V/U$ can only make orbits 
of length $1$ or $q$ on $X$.
\end{proof}
Now given a column in the bottom right quarter corresponding to $K \in \mathcal{R}$ with $\gamma (K) = H$ the marks in the columns corresponding to $H$ and 
$K$ are congruent modulo $q$. The practical significance of Lemmas \ref{ub} and \ref{congruentcolumns} is the following;  Lemma \ref{ub} provides an upper bound 
for each mark in the bottom right quarter. We then utilize Lemma \ref{congruentcolumns} to produce, for each undecided mark in the bottom right quarter, a 
finite range of possible values which the mark might take. It is worth noting that if the upper bound obtained from Lemma \ref{ub} is an integer $< q$ then we immediately 
obtain the correct mark in the bottom right quarter. 

The task now is to attempt to reduce the size of the finite range of values at each undecided position in the bottom right quarter. 
\subsubsection{Transitivity}\label{trans}
Our first tool to reduce the number of possibilities at each position in the bottom right quarter is based on the notion of transitivity. This process provides upper and 
lower bounds for undecided marks in the bottom right quarter of the table of marks of $S$. The procedure, which is described below, is based on the transitivity 
of subgroup inclusion,
\[ U \leq V \ \mathrm{and} \ V \leq K \Rightarrow U \leq K .\]
In terms of conjugacy classes of subgroups this means the following. If $V$ is contained in $p$ conjugates of $K$ then so is $U$. And if $V$ contains $m$ conjugates of $U$ 
then so does $K$.  

At this point in the computation an undecided entry, $\beta_{S/K}(U)$, is represented by a finite range of possible values, one of which is the correct mark. The strategy is to use 
transitivity to reduce the number of values in this range. For clarity we distinguish between the following two situations in Corollary \ref{trans1} and Corollary 
\ref{trans2}.
\begin{cor}\label{trans1}
Let $U \leq V \leq K$. Then 

\begin{enumerate}
 \item any lower bound for $\beta_{S/K}(V)$ is also a lower bound for $\beta_{S/K}(U)$.
\item $ \beta_{S/K}(U) \geq \beta_{S/V}(U)/|K:V| $ .
\end{enumerate}
\end{cor}
\begin{proof}
(i) Follows from Lemma \ref{ub}. (ii)  Follows from the fact that $K$ contains at least as many conjugates of 
$U$ as $V$ does, together with Formula \ref{marksasincidence}.
\end{proof}

\begin{cor}\label{trans2}
Let $V \leq U \leq K$. Then any upper bound for $\beta_{S/K}(V)$ is also an upper bound for $\beta_{S/K}(U)$.
\end{cor}
\begin{proof}
Follows from the fact that $U$ is contained in at least as many conjugates of $K$ as $V$ is, or simply from Lemma \ref{ub}.
\end{proof}

\subsubsection{Dress Congruences}
\label{dresscongruences}
In this section we will describe a refinement of the Dress congruences which enables us to decide the correct entry in many of the positions in the bottom right quarter.
Let $U \leq A$. As before denote $W = N_{S}(U)/U$, and regard $W$ as the union of $B = N_{A}(U)/U$ (its ``blue'' 
elements) and $R = W \setminus B$ (its ``red'' elements). Note that $|B| = \frac{1}{p}|W|$ and that $|R| = (p-1)|B| = \frac{p-1}{p}|W| $. If $X$ is an $S$-set, then 
$Y = \mathrm{Fix}_{X}(U)$ is a $W$-set and by restriction a $B$-set. 

Consider the $S$-set $X = S/K$ for a red subgroup $K$ with $\gamma (K) = H \leq A$. By Proposition \ref{blueonred}, $X$ is equivalent to $A/H$ as an $A$-set. 
It follows that $ \mathrm{Fix}_{S/K}(H)$ is equivalent to $Y = \mathrm{Fix}_{A/H}(H)$ as $B$-sets. We set
\[ o_{W} = \frac{1}{|W|} \sum \limits_{w \in W}^{} \pi_{Y}(w) \]
to be the number of orbits of $W$ on $Y$, and set 
\[ o_{B} = \frac{1}{|B|} \sum \limits_{w \in B}^{} \pi_{Y}(w) \]
to be the number of orbits of $B$ on $Y$. We also set 
\[ o_{R} = \frac{1}{|B|} \sum \limits_{w \in R}^{} \pi_{Y}(w). \]
\begin{prop}
\label{obor}
 With the above notation
\begin{enumerate}
 \item $o_{R} \equiv -o_{B} \pmod{p}$,
\item $o_{R} \leq (p-1)o_{B}$.
\end{enumerate}

\end{prop}
\begin{proof}
 By construction,
\[ p o_{W} = o_{B} + o_{R} \]
and $o_{B} \in \mathbb{Z}$ implies $o_{R} \in \mathbb{Z}$ and 
\[ o_{B} + o_{R} \equiv 0 \pmod{p}. \]
Moreover, $B \leq W$ implies $o_{W} \leq o_{B}$, and thus
\[ o_{R} = p o_{W} - o_{B} \leq p o_{B} - o_{B} = (p - 1)o_{B} \]
as claimed.
\end{proof}
Let $\{ H_{i} \}, i = 1, \ldots, b$ and $\{ K_{j} \}, j = 1, \ldots, r$ be a list of representatives of $\mathcal{B}/S$ and $\mathcal{R}/S$ respectively, and let $X$ be an $S$-set. 
It follows from Theorem \ref{dresscongtheorem} that,
\begin{align}\label{owequalsnuh}
 \sum \limits_{i=1}^{b} n(U, H_{i}) \beta_{X}(H_{i}) + \sum \limits_{j=1}^{r} n(U, K_{j}) \beta_{X}(K_{j}) = c \cdot |W|,  
\end{align}
for $U \leq S$ where $c$ is the number of orbits of $W$ on $Y = \mathrm{Fix}_{X}(U)$, i.e. $c = o_{W}$. Moreover, 
\begin{align}
  \sum \limits_{i=1}^{b} n(U, H_{i}) \beta_{X}(H_{i}) = |B| \cdot o_{B},
\end{align}
and
\begin{align}
 \sum \limits_{j=1}^{r} n(U, K_{j}) \beta_{X}(K_{j}) = |B| \cdot o_{R}.
\end{align}
Since the numbers $o_{B}$ are determined by the marks in the bottom left quarter of the table of marks of $S$, we get the following conditions on the marks in 
the bottom right quarter.
\begin{cor}\label{conginaction}
 Let $K \in \mathcal{R}$ and let $U,o_{B}$ be as above. Then the marks $\beta_{S/K}(K_{j})$ must satisfy,
\[ \frac{1}{|B|} \sum \limits_{j=1}^{r} n(U,K_{j})\beta_{S/K}(K_{j}) \equiv -o_{B}\pmod{p} \]
and
\[ \frac{1}{|B|} \sum \limits_{j=1}^{r} n(U,K_{j})\beta_{S/K}(K_{j}) \leq (p - 1)\cdot o_{B}. \]
\end{cor}

\begin{example}
Table \ref{coeffmat} shows the complete Dress congruence matrix for $S_{5}$. The integer entries in the table represent the numbers 
\[ n(U,H) = \# \{Ua \in N_{S_{5}}(U)/U : \langle U, a \rangle \sim_{S_{5}} H \} \]
where $U$ and $H$ run over a transversal of the conjugacy classes of subgroups of $S$. The final column lists $|W|$ for $W = N_{S_{5}}(U)/U$.

\bigskip

\begin{figure}[H]
\renewcommand{\tabcolsep}{2.5pt}
\begin{tabular}{l|rrrrrrrrr|rrrrrrrrrr|r}

$U$&$ 1$&$ C_{2}$&$ C_{3}$&$ 2^2$&$ C_{5}$&$ S_{3}$&$ D_{10}$&$ A_{4}$&$ A_{5}$&$ C_{2}$&$ C_{4}$&$ 2^2$&$ S_{3}$&$ C_{6}$&$ D_{8}$&$ D_{12}$&$ 5{:}4$&$ S_{4}$&$ S_{5} $&$|W|$\\
\hline
$1 $&$1 $&$ 15 $&$ 20 $&$ × $&$ 24 $&$ × $&$ × $&$ × $&$ × $&$ 10 $&$ 30 $&$  $&$  $&$ 20 $&$ × $&$ × $&$ × $&$ × $&$ ×$&$120$\\
 
$ C_{2} $&$ $&$ 1 $&$ × $&$ 1 $&$ × $&$ × $&$ × $&$ × $&$ × $&$  $&$ 1 $&$ 1 $&$ × $&$ × $&$ × $&$ × $&$ × $&$ × $&$ ×$&$4$\\

$C_{3} $&$  $&$  $&$ 1 $&$ × $&$  $&$1 $&$ × $&$ × $&$ × $&$ × $&$ × $&$  $&$ 1 $&$ 1 $&$ × $&$ × $&$ × $&$ × $&$ ×$&$4$\\

$2^2 $&$ $&$  $&$  $&$ 1 $&$ × $&$ × $&$  $&$ 2 $&$ × $&$ × $&$ × $&$ × $&$ × $&$  $&$ 3 $&$ × $&$ × $&$ × $&$ ×$&$6$\\

$C_{5}$&$ $&$  $&$  $&$  $&$ 1 $&$ × $&$ 1 $&$ × $&$ × $&$ × $&$ × $&$ × $&$ × $&$ × $&$ × $&$ × $&$ 2 $&$ × $&$ ×$&$4$\\

$S_{3}$&$ $&$  $&$  $&$  $&$  $&$ 1 $&$ × $&$ × $&$ × $&$ × $&$ × $&$ × $&$ × $&$ × $&$ × $&$ 1 $&$ × $&$ × $&$ ×$&$2$\\

$D_{10} $&$  $&$  $&$  $&$  $&$  $&$  $&$ 1 $&$ × $&$ × $&$ × $&$ × $&$ × $&$ × $&$ × $&$ × $&$ × $&$ 1 $&$ × $&$ ×$&$2$\\

$A_{4} $&$  $&$  $&$  $&$  $&$  $&$  $&$  $&$ 1 $&$ × $&$ × $&$ × $&$ × $&$ × $&$ × $&$ × $&$ × $&$ × $&$ 1 $&$ ×$&$2$\\

$A_{5} $&$ $&$  $&$  $&$  $&$  $&$  $&$  $&$  $&$ 1 $&$ × $&$ × $&$ × $&$ × $&$ × $&$ × $&$ × $&$ × $&$ × $&$ 1$&$2$\\
\hline

$C_{2}$&$  $&$  $&$  $&$  $&$  $&$  $&$  $&$  $&$  $&$ 1 $&$ × $&$ 3 $&$ × $&$ 2 $&$ × $&$ × $&$ × $&$ × $&$ ×$&$6$\\

$C_{4} $&$ $&$  $&$  $&$  $&$  $&$  $&$  $&$  $&$  $&$  $&$ 1 $&$ × $&$ × $&$ × $&$ 1 $&$ × $&$ × $&$ × $&$ ×$&$2$\\

$2^2$&$ $&$  $&$  $&$  $&$  $&$  $&$  $&$  $&$  $&$  $&$  $&$ 1 $&$ × $&$ × $&$ 1 $&$ × $&$ × $&$ × $&$ ×$&$2$\\

$S_{3}$&$ $&$  $&$  $&$  $&$  $&$  $&$  $&$  $&$  $&$  $&$  $&$  $&$ 1 $&$ × $&$ × $&$1 $&$  $&$ × $&$ ×$&$2$\\

$C_{6}$&$ $&$  $&$  $&$  $&$  $&$  $&$  $&$  $&$  $&$  $&$  $&$  $&$  $&$ 1 $&$ × $&$ 1 $&$ × $&$ × $&$ ×$&$2$\\

$D_{8}$&$ $&$  $&$  $&$  $&$  $&$  $&$  $&$  $&$  $&$  $&$  $&$  $&$  $&$  $&$ 1 $&$ × $&$ × $&$ × $&$ ×$&$1$\\

$D_{12}$&$ $&$  $&$  $&$  $&$  $&$  $&$  $&$  $&$  $&$  $&$  $&$  $&$  $&$  $&$  $&$ 1 $&$ × $&$ × $&$ ×$&$1$\\

$5{:}4$&$ $&$  $&$  $&$  $&$  $&$  $&$  $&$  $&$  $&$  $&$  $&$  $&$  $&$  $&$  $&$  $&$ 1 $&$ × $&$ ×$&$1$\\

$S_{4}$&$ $&$  $&$  $&$  $&$  $&$  $&$  $&$  $&$  $&$  $&$  $&$  $&$  $&$  $&$  $&$  $&$  $&$ 1 $&$ ×$&$1$\\

$S_{5}$&$ $&$ $&$ $&$  $&$ $&$ $&$ $&$  $&$ $&$ $&$ $&$  $&$ $&$ $&$ $&$  $&$ $&$ $&$ 1 $&$1$\\

\end{tabular} 
\caption{Dress Congruence Matrix for $S_{5}$}
\label{coeffmat}
\end{figure}
For example, the congruence corresponding to $U = 1$ is 
\[ y_{1} + 15y_{2} + 20y_{3} + 24y_{5} + 10y_{10} + 30y_{11} + 20y_{14} \equiv 0 \pmod{120}. \]
Each row of the table of marks of $S_{5}$ must satisfy all the congruences. 

To illustrate how Corollary \ref{conginaction} yields conditions on the marks in the bottom right quarter, consider the impression  
\[
 \beta_{S_{5}/D_{12}} = ( 10  ,  2  ,  1  ,  0  ,  0  ,  1  ,  0  ,  0  ,  0  ,  y_{10}  ,  y_{11}  ,  y_{12}  ,  y_{13}  ,  y_{14}  ,  y_{15}  ,  y_{16}  ,  y_{17}  ,  y_{18}  ,  y_{19}) \\
\]
of $S_{5}$ on $S_{5}/D_{12}$.
The marks $\{ y_{1}, \dots, y_{9} \}$ of the blue subgroups are known from Section \ref{sec:bottomrightquarter}. The marks of the red subgroups are represented by $y_{i}$ for 
$i \in \{10, \ldots,19 \}$.
The congruence from $U = C_{2}$ in the top half of Figure \ref{coeffmat} reads
\[ y_{2} + y_{4} + y_{11} + y_{12} \equiv 0 \pmod 4 . \]
Clearly $o_{B} = \frac{1}{2}(y_{2} + y_{4}) = 1$. Moreover, $o_{R} = \frac{1}{2}(y_{11} + y_{12})$. It follows from Corollary \ref{conginaction} that
\begin{enumerate}
 \item $o_{R} \equiv 1 \pmod 2$
\item $o_{R}  \leq 1.$
\end{enumerate}
Hence $o_{R} = 1$ and so $y_{11} + y_{12} = 2$. Lemmas \ref{ub} and \ref{congruentcolumns} yield $y_{11}, y_{12} \in \{0, 2\}$. We 
conclude that either $y_{11} = 0, y_{12} =2$ or $y_{11} = 2, y_{12} = 0$. In this fashion the congruences yield conditions on the marks in the bottom right quarter of the table of marks.

\end{example}
\subsubsection{Explicit Testing of Incidences}
\label{explicit}
If all other approaches fail, one can explicitly count the number of conjugates of $K$ which lie above a subgroup $V$ and compute the mark $\beta_{S/K}(V)$ 
using Proposition \ref{marksasincidence}. 

In order to avoid listing entire conjugacy classes of subgroups, we introduce the following subsets of a conjugacy class of subgroups. For a subgroup $K \leq S$ 
and an element $t \in S$ denote 
\[ X(K,t) = \{ K' \in [K]_{S} : t \in K' \}. \]

\begin{lemma}
 Let $V \leq S$ and $t \in S$. Then
\[ \{ K' \in [K]_{S} : V \leq K' \} = \{ K' \in X(K, t) : V \leq K' \} \]
\end{lemma}
 \begin{proof}
  By definition $X(K, t)$ is precisely the subset of $[K]_{S}$ consisting of those conjugates $K'$ of $K$ which contain the element $ t \in V$. Thus $K' \geq V$ 
implies $K' \in X(K,t)$.
 \end{proof}
\noindent In particular if $V$ is a red subgroup and $t \in V \setminus A$ then 
\[ \beta_{S/K}(V) = |N_{S}(K):K| \cdot \# \{ K' \in X(K,t) : V \leq K' \}. \]
Such a set $X(K,t)$ can be computed efficiently, using Proposition \ref{pairs2}, as follows.
\begin{prop}
Let $K \leq S$ and $t \in S$. Then
\begin{enumerate}
 \item the centralizer $C = C_{S}(t)$ acts on $X(K,t)$ by conjugation;
\item the normalizer $N = N_{S}(K)$ acts on $T = K \cap [t]_{S}$ by conjugation;
\item the map $\xi : X/C \rightarrow T/N$ given by
\[ \xi([K^s]_{C}) = [t^{s^{-1}}]_{N} \]
is a well defined bijection.
\end{enumerate}

\end{prop}
 \begin{proof}
(i) and (ii) are obvious. (iii) If $Z = \{ (K', t') \in [K]_{S} \times [t]_{S} : t' \in K' \}$ then $Z$ is $S$-invariant, $X(k,t) = Zt$ and the claim follows with Proposition \ref{pairs2}.
 \end{proof}
\noindent This result allows us to compute the set 
\[ X(K,t) = \coprod \limits_{[a]_{N} \in K/N, a^s = t}^{} [K^s]_{C}\]
systematically as a disjoint union of $C$-orbits of conjugates of $K$, by first computing the conjugacy classes of elements of $K$, partitioning them 
into $N$-orbits, and selecting those consisting of conjugates of $t$. For each such $N$-orbit $[a]_{N}$ one finds a conjugating element $s \in S$ with $a^s = t$ and then 
computes the $C$-orbit of the conjugate $K^s$.

\section{Computation}\label{sec:computation}
Propositions \ref{blueonblue} and \ref{blueonred} enable us to determine the marks in the top left and bottom left quarters respectively. The bounds 
described in Section \ref{bounds} yield a partially complete bottom right quarter, where, if a mark is undecided, it is represented by a finite range of values. 
We work our way down through the table of marks completing each row before we move on to the next one. We apply the congruences and the transitivity tests 
until the row is completed or no new mark is obtained. If there are still undecided marks we use the explicit incidence test from  Section \ref{explicit} with a single 
$t$ to compute as many marks as possible. Then we apply the congruences and transitivity tests again. If there are still undecided marks we run the incidence 
test again with a different $t$ and repeat the process until the row is complete. The entire process is summarized in Algorithm \ref{tomalg}.
\begin{algorithm}[H]                 
\caption{ \texttt{TableOfMarksByCyclicExtension()}} \label{tomalg}         
\begin{algorithmic}[H]
\State \textbf{Input} Subgroup pattern $( \mathrm{Sub}(A)/A, \mathrm{M}(A))$ of $A$.
\State \textbf{Output} Subgroup pattern of $S$.
\State Compute $\mathrm{Sub}(S)/S$ as \texttt{SubgroupsByCyclicExtension(}$\mathrm{Sub}(A)/A$\texttt{)}.
\State Use Proposition \ref{blueonblue} to compute top left quarter of $\mathrm{M}(S)$.
\State Use Proposition \ref{blueonred} to compute bottom left quarter of $\mathrm{M}(S)$. 
\For{ each row in bottom right of $\mathrm{M}(S)$}
\State Implement bounds from Subsection \ref{bounds}.
\While{row is incomplete}
\State Apply congruences (\ref{dresscongruences}) and 
\State transitivity (\ref{trans}) until no more new marks are found.
\If{ row still contains undecided marks}
\State Compute some marks explicitly (\ref{explicit}).
\EndIf
\EndWhile
\EndFor
\State \textbf{return} $(\mathrm{Sub}(S)/S, \mathrm{M}(S))$.
\end{algorithmic}
\end{algorithm}

This algorithm completes the task of computing the subgroup pattern of $S$ from that of $A$. Some of the results obtained by a \texttt{GAP} implementation of this 
algorithm are listed in Section \ref{sec:results}.

\subsection{Computing the Table of Marks of a Solvable Group}
In Section \ref{solvablesubs} we described a new algorithm to compute the conjugacy classes of subgroups of a solvable group $G$. In the same spirit we have developed an
 algorithm to compute the table of marks of a solvable group $G$ based on the procedures described in the preceding sections. The strategy is the same as in 
Section \ref{solvablesubs}. We take as input a solvable group $G$, and work our way up through the composition series of $G$ starting with the table of marks of the 
trivial group, computing the table of marks of each group in the series in turn until we obtain the table of marks of $G$ itself.

\begin{algorithm}[H]                
\caption{ \texttt{TableOfMarksSolvableGroup()}}
    \label{tomsolvable}       
\begin{algorithmic}[H]
\State \textbf{Input} A solvable group $G$.
\State \textbf{Output} Subgroup pattern $(\mathrm{Sub}(G)/G, \mathrm{M}(G))$ of $G$.
\State Compute a composition series $1 = G_{0} \trianglelefteq G_{1} \trianglelefteq \ldots \trianglelefteq G_{n} = G$
\State Set $P_{0} \gets ( \mathrm{Sub}(1)/1, \mathrm{M}(1))$.
\For{ $i \in \{1, \ldots, n \}$ }
\State $P_{i} \gets$ \texttt{TableOfMarksByCyclicExtension(}$P_{i-1}$\texttt{)}.
\EndFor
\State \textbf{return} $P_{n}$.
\end{algorithmic}
\end{algorithm}
\begin{example}
Recall the example of $\mathrm{GL}_{2}(3)$ from Section \ref{solvablesubs}, and its associated composition series
\[1  \vartriangleleft 2 \vartriangleleft 4 \vartriangleleft Q_{8} \vartriangleleft \textrm{SL}_{2} (3) \vartriangleleft \mathrm{GL}_{2} (3) \]
 In this example we apply Algorithm \ref{tomsolvable} starting with the 
table of marks of the trivial group to obtain the table of marks of $\mathrm{GL}_{2}(3)$.
\begin{figure}[H]
 \[ \left( \begin{array}{c}
1  
\end{array} \right)
\overset{p = 2}{\longrightarrow}
\left( \begin{array}{c|c}
2\\
\hline
1& 1
\end{array} \right)
\overset{p = 2}{\longrightarrow}
\left( \begin{array}{cc|c}
 4& \\
 2 &2\\
\hline
 1 &1& 1
\end{array} \right)
\overset{p = 2}{\longrightarrow}
\left( \begin{array}{ccc|ccc}
 8& & \\
 4 &4 &\\
 2 &2 &2\\
\hline
 2 &2 &. &2\\
 2 &2 &.& .& 2\\
 1 &1 &1 &1 &1 &1

\end{array} \right) \overset{p = 3}{\longrightarrow}
\]

\[\arraycolsep 3.7pt \left( \begin{array}{cccc|ccc}

 24& & &\\
 12 &12& &\\
 6  &6& 2 &\\
 3  &3& 3& 3\\
\hline
 8 & .& .& .& 2\\
 4 & 4 &.& .& 1& 1\\
 1 & 1& 1& 1& 1& 1& 1
\end{array} \right) 
\overset{p = 2}{\longrightarrow}              
\left( \begin{array}{ccccccc|ccccccccc}
 48& & & & & & \\
 24& 24& & & & & \\
 16&  .& 4& & & & \\
 12& 12& .& 4& & & \\
 8 & 8& 2& .& 2& & \\
 6 & 6 &.& 6& .& 6& \\
 2 & 2& 2& 2& 2& 2& 2 \\
\hline
 24&  .& .& .& .& .& . &2\\
 12& 12& .& .& .& .& .& 2& 2\\
 8&  .& 2& .& .& .& .& 2& .& 2\\
 8&  .& 2& .& .& .& .& 2& .& .& 2\\
 6&  6& .& 2& .& .& .& 2& 2& .& .& 2\\
 6&  6& .& 2& .& .& .& .& .& .& .& .& 2\\
 4&  4& 1& .& 1& .& .& 2& 2& 1& 1& .& .& 1\\
 3&  3& .& 3& .& 3& .& 1& 1& .& .& 1& 1& .& 1\\
 1&  1& 1& 1& 1& 1& 1& 1& 1& 1& 1& 1& 1& 1& 1& 1
\end{array} \right)
\]
\caption{Table of Marks of $\mathrm{GL}_{2}(3)$}
\end{figure}
\end{example}

\subsection{Results and Statistics}\label{sec:results}
The methods described in this article have been used to extend the GAP table of marks library Tomlib. Tables \ref{results} and \ref{more-results} list 
some of the groups to which these methods have been applied together with running times for the computations. Table \ref{results} contains two extra columns labeled 
$\# X(K,t)$ and $\mathrm{max}|X(K,t)|$ where $\# X(K,t)$ records the number of times a mark is computed explicitly based on Section \ref{explicit}, and $\mathrm{max}|X(K,t)|$ 
records the length of the largest orbit which is computed for such a calculation. The computations were carried out on an 
Apple MacBook Pro with an Intel Core 2 Duo CPU T7500 @ 2.20GHz with 2 gigabytes of RAM.
\bigskip
\begin{table}[H]
\begin{center}
\begin{tabular}{|l|l|r|r|r|r|r|}
\hline
$A$ & $S$ & $|\mathrm{Sub}(A)/A|$ & $|\mathrm{Sub}(S)/S|$ &$\# X(K,t)$&$\mathrm{max}\left|X(K,t)\right|$& $\mathrm{Time}$\\
\hline
$A_5$ & $S_5$ & $9$ & $19$ &$0$ &$0$ & $1$s\\
$A_6$ & $S_6$ & $22$ & $56$ &$2$ &$4$ & $2$s\\
$A_7$ & $S_7$ & $40$ & $96$ &$3$ &$20$ & $3$s\\
$A_8$ & $S_8$ & $137$ & $296$ &$26$ &$60$ & $20$s\\
$A_9$ & $S_9$ & $223$ & $554$ &$82$ &$140$ & $50$s\\
$A_{10}$ & $S_{10}$ & $430$ & $1593$ &$381$ &$384$ & $6$m\\
$A_{11}$ & $S_{11}$ & $788$ & $3094$ &$912$ &$960$ & $20$m\\
$A_{12}$ & $S_{12}$ & $2537$ & $10723$ & $6161$ & $3240$ & $7$h\\
$A_{13}$ & $S_{13}$ & $4558$ & $20832$ & $12316$ & $15120$ & $43$h\\
\hline
\end{tabular}
\end{center}
\caption{Results for Symmetric Groups}
\label{results}
\end{table}
\bigskip
\begin{table}[H]
\begin{center}
\begin{tabular}{|l|l|r|r|r|}
\hline
$A$ & $S$ & $|\mathrm{Sub}(A)/A|$ & $|\mathrm{Sub}(S)/S|$ & $\mathrm{Time}$\\
\hline
$\mathrm{He}$ & $\mathrm{He}.2$ & $1698$ & $1930$ & $231$m\\
$\mathrm{HS}$ & $\mathrm{HS}.2$ & $589$ & $2057$ & $35$m\\
$Sz(8)$ & $Sz(8).3$ & $22$ & $39$ & $3s$\\
$^{2}F_{4}(2)'$ & $^{2}F_{4}(2)$ & $434$ & $849$ & $48$m\\
$L_{2}(32)$ & $L_{2}(32).5$ & $24$ & $30$ & $4$s\\
\hline
\end{tabular}
\end{center}
\caption{More Results}
\label{more-results}
\end{table}
A \textsf{GAP} implementation of the algorithms is available on request from the authors. 

\bigskip \noindent{\bf Acknowledgment}: 
Much of the work in this article is based on the first authors PhD thesis (see \cite{thesis}). This research  was supported by  Science Foundation
Ireland (07/RFP/MATF466).


\bibliographystyle{amsplain}
\bibliography{References}

\def\cprime{$'$}
\providecommand{\bysame}{\leavevmode\hbox to3em{\hrulefill}\thinspace}
\providecommand{\MR}{\relax\ifhmode\unskip\space\fi MR }
\providecommand{\MRhref}[2]{%
  \href{http://www.ams.org/mathscinet-getitem?mr=#1}{#2}
}
\providecommand{\href}[2]{#2}
\begin{thebibliography}{10}

\bibitem{bouc}
Serge Bouc, \emph{Burnside rings}, Handbook of algebra, {V}ol. 2,
  North-Holland, Amsterdam, 2000, pp.~739--804. \MR{1759611 (2001m:19001)}

\bibitem{breuer}
Thomas Breuer and G{\"o}tz Pfeiffer, \emph{Finding possible permutation
  characters}, J. Symbolic Comput. \textbf{26} (1998), no.~3, 343--354.
  \MR{1633876 (99e:20005)}

\bibitem{geometriesandgroups}
Francis Buekenhout, \emph{Diagrams for geometries and groups}, J. Combin.
  Theory Ser. A \textbf{27} (1979), no.~2, 121--151. \MR{542524 (83f:51003)}

\bibitem{structurem12}
Francis Buekenhout and Sarah Rees, \emph{The subgroup structure of the
  {M}athieu group {$M_{12}$}}, Math. Comp. \textbf{50} (1988), no.~182,
  595--605. \MR{929556 (88m:20024)}

\bibitem{burn}
W.~Burnside, \emph{Theory of groups of finite order}, Dover Publications Inc.,
  New York, 1955, 2d ed. \MR{0069818 (16,1086c)}

\bibitem{dress}
Andreas Dress, \emph{A characterisation of solvable groups}, Math. Z.
  \textbf{110} (1969), 213--217. \MR{0248239 (40 \#1491)}

\bibitem{GAP4}
The GAP~Group, \emph{{GAP -- Groups, Algorithms, and Programming, Version
  4.4.12}}, 2008, \verb+http://www.gap-system.org+.

\bibitem{hulpke}
Alexander Hulpke, \emph{Computing subgroups invariant under a set of
  automorphisms}, J. Symbolic Comput. \textbf{27} (1999), no.~4, 415--427.
  \MR{1681348 (2000a:20001)}

\bibitem{thesis}
Liam Naughton, \emph{Computing the table of marks of a finite group}, Ph.D.
  thesis.

\bibitem{neumann}
Peter~M. Neumann, \emph{A lemma that is not {B}urnside's}, Math. Sci.
  \textbf{4} (1979), no.~2, 133--141. \MR{562002 (81g:01012)}

\bibitem{pfe}
G{\"o}tz Pfeiffer, \emph{The subgroups of {$M_{24}$}, or how to compute the
  table of marks of a finite group}, Experiment. Math. \textbf{6} (1997),
  no.~3, 247--270. \MR{1481593 (98h:20032)}

\bibitem{tomlib}
\emph{{Tomlib, Version 1.2.1 }},  (2011),
  \verb+http://schmidt.nuigalway.ie/tomlib+.

\end{thebibliography}
\end{document}